%% file: main_rev_02.tex
\documentclass[10pt,journal]{IEEEtran}
\usepackage{graphicx}          
\usepackage{amsmath}           
\usepackage{amssymb}
\usepackage{algorithm}
\usepackage{algpseudocode}
\usepackage{psfrag}
\usepackage[latin1]{inputenc}  
\usepackage{hyperref}

\usepackage{amsthm}

\theoremstyle{definition}

\theoremstyle{remark}

\theoremstyle{plain}

\newcommand{\mbf}[1]{\mathbf{#1}}

\begin{document}

\title{Bayesian Estimation with Distance Bounds}
\author{Dave Zachariah, Isaac Skog, Magnus Jansson, and Peter Händel\thanks{The authors are with the ACCESS Linnaeus Centre, KTH Royal Institute of
Technology, Stockholm, Sweden. E-mail:
$\{$dave.zachariah,isaac.skog,magnus.jansson,peter.handel$\}$@ee.kth.se.}}


\maketitle

\begin{abstract}
We consider the problem of estimating a random state vector when there is information about the maximum distances between its subvectors. The estimation problem is posed in a Bayesian framework in which the minimum mean square error (MMSE) estimate of the state is given by the conditional mean. Since finding the conditional mean requires multidimensional integration, an approximate MMSE estimator is proposed. The performance of the proposed estimator is evaluated in a positioning problem. Finally, the application of the estimator in inequality constrained recursive filtering is illustrated by applying the estimator to a dead-reckoning problem. The MSE of the estimator is compared with two related posterior Cramér-Rao bounds.
\end{abstract}

\input{intro}
\input{problem}
\input{solution}
\input{experiment}
\input{conclusions}

\IEEEtriggeratref{100}
\bibliography{refs_bayes}
\bibliographystyle{ieeetr}

\end{document}

%% file: intro.tex
\section{Introduction}
Consider the problem of estimating a vector $\mbf{x}
\in \mathbb{R}^d$ with a known prior distribution. In the standard Bayesian setup one tries to estimate $\mbf{x}$ from a correlated observation $\mbf{y}$. In this letter we deviate from the standard problem and instead seek the estimate of $\mbf{x}$ relying on the side information  $\|\mbf{x}_i - \mbf{x}_j\| \leq \gamma$, where $\mbf{x}_i$ and $\mbf{x}_j$ are subvectors of $\mbf{x}$, and $\gamma$ is given.

This setup has a range of applications in positioning
and localization. As an example, consider the problem of
estimating jointly the positions $\mbf{p}_1\in \mathbb{R}^3$ and $\mbf{p}_2 \in \mathbb{R}^3$ and the velocities $\mbf{v}_1\in \mathbb{R}^3$ and $\mbf{v}_2\in \mathbb{R}^3$ of two points on a
human body. In that case, the state is $\mbf{x} = [\mbf{p}^\top_1 \quad
\mbf{p}^\top_2 \quad \mbf{v}^\top_1 \quad \mbf{v}^\top_2  ]^\top$,
and since there is an upper limit $\gamma$ on how far apart body
parts can be, we have the side information that
$\|\mbf{p}_1-\mbf{p}_2\| \leq \gamma$.
Similarly, consider the problem of
estimating the positions of $N$ nodes, $\mbf{x} = [\mbf{x}^\top_1
\dots \mbf{x}^\top_N  ]^\top$, when bounds, $\|\mbf{x}_i-\mbf{x}_j\| \leq \gamma_{ij}$, on the distance between pairs of nodes are given. Further applications are possible for sensor fusion of systems subject to non-rigid constraints.

In the literature, there are two main approaches to tackle the problem of estimation with nonlinear inequality constraints \cite{Simon2010}. The first approach uses the side information by passing $\mbf{x}$ through a nonlinear function, such that the output always fulfills the constraint \cite{Julier&LaViola2007, Lan&Li2011}. The second approach treats the problem probabilistically, and the side information is used to form a conditional probability density function (pdf) of $\mbf{x}$ by truncating and renormalizing the prior pdf; the support of the conditional pdf is a region where the constraint is inactive \cite{Simon&Simon2010, StrakaEtAl2012}. However, when using distance bounds as side information, the support of the conditional pdf is infinite. In scenarios where the dispersion of the pdf increases without bound, this may lead to numerical problems when computing the moments of the pdf using e.g., the Monte Carlo methods suggested in \cite{StrakaEtAl2012}.

In this letter, we circumvent the problems associated with the infinite support pdf, by reformulating the problem via a linear transformation. This also reduces the dimensionality of the original problem by a factor two. Based on the new problem formulation an efficient method for computing an approximate minimum mean square error (MMSE) estimator of $\mbf{x}$ given bounds on the distances between its subvectors is presented. The method is validated using simulations.

\emph{Notation:} $\mbf{A} \oplus \mbf{B}$ denotes the direct sum
between matrices $\mbf{A}$ and $\mbf{B}$. $\mbf{A}^{1/2}$ is the
lower-triangular Cholesky factorization of a positive definite matrix
$\mbf{A}$. $\text{tr} \{ \cdot \}$ denotes the trace operator, $[\mbf{A}]_i$ is the $i$th column of $\mbf{A}$, and $\mbf{A}^{-\top}$ denotes the transpose of the inverse matrix $\mbf{A}^{-1}$, i.e., $\mbf{A}^{-\top} = (\mbf{A}^{-1})^\top$. Further, $\mbf{I}_n$, $\mbf{1}_n$, and $\mbf{0}_{n}$ denote the identity matrix, a vector of ones and zeros, respectively.


%% file: problem.tex
\section{Problem formulation}
Let the state be defined as $\mbf{x} \triangleq [\mbf{x}^\top_1 \quad
\mbf{x}^\top_2 \quad \mbf{x}^\top_a ]^\top \in \mathbb{R}^d$, where
$\mbf{x}_1, \mbf{x}_2 \in \mathbb{R}^n$ are the two subvectors related to the side information, $\mbf{x}_a \in \mathbb{R}^m$
is a subvector holding the auxiliary states, and the state dimension $d = 2n + m$. Further, we assume that the prior pdf $p(\mbf{x})$ of the state is Gaussian with known mean $\mbf{m}_x$ and covariance $\mbf{C}_x$, i.e., $\mbf{x}\sim\mathcal{N}(\mbf{m}_x,\mbf{C}_x)$.

Now, if we in addition to $\mbf{m}_x$ and $\mbf{C}_{x}$ are provided with the side information $c$ that tells us about the maximum distance between the two subvectors, i.e., we have the constraint that $\| \mbf{Lx} \|\leq \gamma$, where $\mbf{L} \triangleq \begin{bmatrix}\mbf{I}_n & -\mbf{I}_n & \mbf{0}_{n\times m} \end{bmatrix}$. Then, we would like to compute the MMSE estimate of $\mbf{x}$ given $c$. The estimator and its error covariance matrix are given by the mean $\mbf{m}_{x|c}$ and covariance matrix $\mbf{C}_{x|c}$ of
the conditional pdf $p(\mbf{x}|c)$ \cite{Kay1993}.
Since the computation of these moments requires multidimensional numerical integration, our aim is to find a computationally efficient way to calculate $\mbf{m}_{x|c}$ and $\mbf{C}_{x|c}$ approximately. The resulting approximations will be denoted as $\hat{\mbf{x}}_{|c}$ and $\mbf{C}_{\hat{{x}}}$.

%% file: solution.tex
\section{Proposed solution}
The proposed method to calculate $\mbf{m}_{x|c}$ and $\mbf{C}_{x|c}$ consists two steps. In the first step, a state transformation is performed, that maps the infinite integration area into a  finite area and reduces the dimension of the integrals. Thus the transformation simplifies the calculation of the conditional mean and covariance. In the second step, the integrals are approximated using a deterministic sampling approach.

\subsection{Variable transformation}
Define a new state vector, given by the invertible linear
transformation $\mbf{z} \triangleq \mbf{Tx} \in \mathbb{R}^d$,
such that it can be divided into the subvectors $\mbf{z}_1 = \mbf{x}_1
- \mbf{x}_2 \in \mathbb{R}^n$ and $\mbf{z}_2 = [ (\mbf{x}_1 +
\mbf{x}_2)^\top \quad \mbf{x}^\top_a ]^\top \in \mathbb{R}^{n+m}$.
Hence, it can be verified that
\begin{equation*}
\mbf{T} \triangleq \begin{bmatrix} \mbf{I}_n & -\mbf{I}_n  \\
\mbf{I}_n & \mbf{I}_n  \\
\end{bmatrix}
\oplus\mbf{I}_m,\quad\mbf{T}^{-1} = \frac{1}{2}\left(
\begin{bmatrix}
\mbf{I}_n & \mbf{I}_n \\
-\mbf{I}_n & \mbf{I}_n
\end{bmatrix}\oplus2\mbf{I}_m\right),
\end{equation*}
and $\mbf{z} \sim \mathcal{N}(\mbf{m}_z, \mbf{C}_z)$. Here
$\mbf{m}_z = \mbf{Tm}_x$ and $\mbf{C}_z =
\mbf{T}\mbf{C}_x\mbf{T}^\top$. In terms
of the new state vector, the constraint $c$ can be expressed as $ \| \mbf{L T}^{-1}
\mbf{z}\| =  \| \mbf{z}_1 \| \leq \gamma$, and the constraint thus
operates only on the subvector $\mbf{z}_1$. As will been shown next,
this important property simplifies the
calculations of $\mbf{m}_{x|c}$ and $\mbf{C}_{x|c}$.
Since the means and covariances of the new and original state are related via $\mbf{m}_{x} = \mbf{T}^{-1} \mbf{m}_{z}$ and $\mbf{C}_{x} = \mbf{T}^{-1} \mbf{C}_{z} \mbf{T}^{-\top}$, we can recast the original problem of calculating $\mbf{m}_{x|c}$ and $\mbf{C}_{x|c}$, into computing $\mbf{m}_{z|c}$ and $\mbf{C}_{z|c}$.

To compute the moments, note that $p(\mbf{z}_2| \mbf{z}_1)$ is Gaussian with mean  $\mbf{m}_{z_2|z_1}$ and covariance $\mbf{C}_{z_2|z_1}$ given by the affine mapping
    \begin{equation}
    \begin{split}
        \mbf{m}_{z_2|z_1} &= \mbf{u}_{z_2} + \mbf{A} \mbf{z}_1 \\
        \mbf{C}_{z_2|z_1} &= \mbf{C}_{z_2} - \mbf{A}\mbf{C}^\top_{z_2 z_1},
    \end{split}
    \end{equation}
where $\mbf{A} \triangleq \mbf{C}_{z_2 z_1} \mbf{C}^{-1}_{z_1}$ and $\mbf{u}_{z_2} \triangleq \mbf{m}_{z_2} - \mbf{A} \mbf{m}_{z_1}$.

Next, since $\mbf{z}_1$ is \emph{given} the norm $\| \mbf{z}_1 \|$ provides no additional information and $p(\mbf{z}_2| \mbf{z}_1, \| \mbf{z}_1 \|)$ is identical to $p(\mbf{z}_2| \mbf{z}_1)$. Similarly, given $\mbf{z}_1$ a bound on the norm $\| \mbf{z}_1 \| \leq \gamma$ yields no additional information. Hence $p(\mbf{z}_2|\mbf{z}_1, c) = p(\mbf{z}_2 | \mbf{z}_1)$ for all valid $\mbf{z}_1$, i.e., $\forall\mbf{z}_1\in\{\mbf{z}_1\in \mathbb{R}^n : \|\mbf{z}_1\|\leq \gamma\}$. The conditional mean $\mbf{m}_{z|c}$ and covariance $\mbf{C}_{z|c}$ can thus be calculated as follows.

Let the conditional mean be partitioned as $\mbf{m}_{z|c} =
[\mbf{m}^\top_{z_1|c} \quad \mbf{m}^\top_{z_2|c}]^\top$. The conditional mean of $\mbf{z}_1$ is given by
\begin{equation}
\begin{split}
\mbf{m}_{z_1|c} 
&= \int_{\mbf{z}_1} \mbf{z}_1 p(\mbf{z}_1|c) d\mbf{z}_1, 
\end{split}
\label{eq:conditionalmean}
\end{equation}
and the conditional mean of $\mbf{z}_2$ by
\begin{equation}
\begin{split}
\mbf{m}_{z_2|c} 
&= \int_{\mbf{z}_1} \left[  \int_{\mbf{z}_2} \mbf{z}_2 p(\mbf{z}_2|\mbf{z}_1) d\mbf{z}_2 \right] p(\mbf{z}_1|c) d\mbf{z}_1 \\
&= \int_{\mbf{z}_1} \left( \mbf{u}_{z_2} + \mbf{A} \mbf{z}_1 \right) p(\mbf{z}_1|c) d\mbf{z}_1 \\
&= \mbf{u}_{z_2}  +  \mbf{A} \mbf{m}_{z_1|c}.
\end{split}
\label{eq:affinem2}
\end{equation}
The conditional covariance matrix can be written as
\begin{equation}
\begin{split}
\mbf{C}_{z|c} &= \mbf{P}_{z|c} - \mbf{m}_{z|c} \mbf{m}_{z|c}^\top \\
&=
\begin{bmatrix}
\mbf{P}_{z_1|c}     & \mbf{P}_{z_1 z_2|c} \\
\mbf{P}^\top_{z_1 z_2|c} & \mbf{P}_{z_2|c}
\end{bmatrix}
- \mbf{m}_{z|c} \mbf{m}_{z|c}^\top,
\end{split}
\label{eq:covariance}
\end{equation}
where
\begin{equation}
\begin{split}
\mbf{P}_{z_1|c} 
&= \int_{\mbf{z}_1} \mbf{z}_1 \mbf{z}^\top_1 p( \mbf{z}_1 | c) d\mbf{z}_1,
\end{split}
\label{eq:conditionalcovariance}
\end{equation}

\begin{equation}
\begin{split}
\mbf{P}_{z_1z_2|c} 
&= \int_{\mbf{z}_1} \mbf{z}_1  \left[ \int_{\mbf{z}_2} \mbf{z}^\top_2
  p(\mbf{z}_2|\mbf{z}_1) d\mbf{z}_2 \right] p(\mbf{z}_1 | c) d\mbf{z}_1 \\
&= \mbf{m}_{z_1|c} \mbf{u}^\top_{z_2} + \mbf{P}_{z_1|c}
\mbf{A}^\top,
\end{split}
\label{eq:affinePz1z2}
\end{equation}
and
\begin{equation}
\begin{split}
\mbf{P}_{z_2|c} 
&= \int_{\mbf{z}_1} \left[  \int_{\mbf{z}_2} \mbf{z}_2 \mbf{z}^\top_2 p(\mbf{z}_2|\mbf{z}_1) d\mbf{z}_2 \right] p(\mbf{z}_1|c) d\mbf{z}_1 \\
&= \int_{\mbf{z}_1} \left( \mbf{C}_{z_2|z_1} + \mbf{m}_{z_2|z_1}\mbf{m}^\top_{z_2|z_1} \right) p(\mbf{z}_1|c) d\mbf{z}_1 \\
&= \mbf{C}_{z_2} - \mbf{A}\mbf{C}^\top_{z_2 z_1} + \mbf{u}_{z_2}\mbf{u}^\top_{z_2} \\
&+\int_{\mbf{z}_1} \left( \mbf{u}_{z_2}\mbf{z}^\top_1 \mbf{A}^\top + \mbf{A}\mbf{z}_{1}\mbf{u}^\top_{z_2} + \mbf{A}\mbf{z}_1 \mbf{z}^\top_1 \mbf{A}^\top\right) p(\mbf{z}_1|c) d\mbf{z}_1 \\
&= \mbf{C}_{z_2} - \mbf{A}\mbf{C}^\top_{z_2 z_1} +
\mbf{u}_{z_2}\mbf{u}^\top_{z_2} \\&+
\mbf{u}_{z_2}\mbf{m}^\top_{z_1|c} \mbf{A}^\top +
\mbf{A}\mbf{m}_{z_1|c}\mbf{u}^\top_{z_2} + \mbf{A}\mbf{P}_{z_1|c}
\mbf{A}^\top.
\end{split}
\label{eq:affinePz2}
\end{equation}
Hence, we need only solve the integrals \eqref{eq:conditionalmean} and
\eqref{eq:conditionalcovariance} for the first and second order moments of $p(\mbf{z}_1|c)$. The remaining parts of $\mbf{m}_{z|c}$ and $\mbf{C}_{z|c}$ can be found by a series of affine transformations. An important computational advantage of the reparameterization, $\mbf{z} = \mbf{Tx}$, is that the support of $p(\mbf{z}_1|c)$ is finite unlike the support of $p(\mbf{x}|c)$. Further, a reduction of the dimensionality by a factor two compared to the original problem is achieved. When $n=1$, the integrals \eqref{eq:conditionalmean} and \eqref{eq:conditionalcovariance} are given by the mean
and variance of a truncated Gaussian distribution, which both have closed-form expressions \cite{JohnsonEtAl1994}.

\subsection{Approximating the conditional mean and covariance}
When $n>1$, no apparent closed-form expressions of integrals
\eqref{eq:conditionalmean} and \eqref{eq:conditionalcovariance}
exist. To avoid solving the integrals through computationally
complex numerical integrations, we will approximate them by the convex combinations
\begin{equation}
\begin{split}
\mbf{m}_{z_1|c} \simeq \sum^{2n}_{i=0} w_i \mbf{z}^{(i)}_1 \quad\mbox{and}\quad
\mbf{P}_{z_1|c} \simeq \sum^{2n}_{i=0} w_i \mbf{z}^{(i)}_1\!
(\mbf{z}^{(i)}_1)^\top
\end{split},
\label{eq:approximation}
\end{equation}
where the $i$th sample point $\mbf{z}^{(i)}_1$ and weight $w_i$, described next, are chosen so that the following properties hold true. When a fraction $\alpha$ of the probability mass of $p(\mbf{z}_1)$
is within the boundary of the constraint, the approximated moments are unchanged. Otherwise, the sample points adapt to the constraint, which ensures that the approximated mean falls within its convex boundary and that the dispersion is reduced accordingly. Thus, the design parameter $\alpha$ affects when the side information becomes effective.


The proposed deterministic sampling technique is as follows. First, as in the
sigma-point transformation \cite{Candy2009}, $2n+1$ sample points
from $p(\mbf{z}_1)$ are generated deterministically by
\begin{equation}
\mbf{s}^{(i)} = \begin{cases} \mbf{m}_{z_1} & i = 0\\
\mbf{m}_{z_1} + \eta^{1/2}_{\alpha} [ \mbf{C}^{1/2}_{z_1} ]_{i} & i = 1, \dots, n \\
\mbf{m}_{z_1} - \eta^{1/2}_{\alpha} [ \mbf{C}^{1/2}_{z_1} ]_{i-n} & i = n+1, \dots, 2n
\end{cases}.
\end{equation}
Here $\eta_\alpha$ is the value that fulfills $\Pr \{ \eta \leq \eta_\alpha \} =\alpha$, where $\eta =
(\mbf{z}_1 - \mbf{m}_{z_1})^{\top} \mbf{C}^{-1}_{z_1}(\mbf{z}_1 -
\mbf{m}_{z_1})$. That is, $\eta_\alpha$ is set by the confidence ellipse at level
$\alpha$ and can be calculated from the inverse of the cumulative
distribution function (cdf) of $\eta \sim \chi^2_n$.

Then, a new set of sample points $\{\mbf{z}^{(i)}_1\}_{i=0}^{2n}$ are generated by orthogonally projecting the sample points $\mbf{s}^{(i)}$ that violate the constraint $c$ onto its spherical boundary. That is
\begin{equation}
\mbf{z}^{(i)}_1 = \begin{cases}
\mbf{s}^{(i)} & \text{if }  \| \mbf{s}^{(i)} \| \leq \gamma \\
\frac{\gamma}{\| \mbf{s}^i \|} \mbf{s}^{(i)} & \text{otherwise}
\end{cases},\quad i=0,\ldots,2n.
\end{equation}
The constraint-violating points are, in a fashion similar to that in~\cite{Lan&Li2011}, resampled at the boundary. The weights in \eqref{eq:approximation} are set as
\begin{equation}
w_i = \begin{cases} 1 - \frac{n}{\eta_\alpha} & i = 0 \\
\frac{1}{2\eta_\alpha} & i = 1, \ldots, 2n
\end{cases}.
\end{equation}
This choice yields the properties for \eqref{eq:approximation} described earlier.\footnote{An $\alpha$ close to one yields a confidence ellipse that captures a larger probability mass and thus adapts smoother to the constraint. However, setting $\alpha$ too large results in sample points with too low weights to approximate the truncated pdf.} Once the moments are computed, $\mbf{m}_{x|c} = \mbf{T}^{-1}
\mbf{m}_{z|c}$ and $\mbf{C}_{x|c} = \mbf{T}^{-1} \mbf{C}_{z|c}
\mbf{T}^{-\top}$ are obtained using the affine transformations
\eqref{eq:affinem2}, \eqref{eq:affinePz1z2}, and
\eqref{eq:affinePz2} along with \eqref{eq:covariance}. This provides the estimator
$\hat{\mbf{x}}_{|c}$ and approximated error covariance matrix
$\mbf{C}_{\hat{x}}$.

\subsection{Illustration of how the estimator works}
\begin{figure}[!t]
  \begin{center}
  \psfrag{mx1}[c][c][0.7][0]{$\mbf{m}_{x_1}$}
  \psfrag{mx2}[c][c][0.7][0]{$\mbf{m}_{x_2}$}
  \psfrag{m1c}[c][c][0.7][0]{$\mbf{m}_{x_1|c}$}
  \psfrag{m2c}[c][c][0.7][0]{$\mbf{m}_{x_2|c}$}
    \psfrag{xaxis}[c][c][0.7][0]{$x$ [m]}
    \psfrag{yaxis}[c][c][0.7][0]{$y$ [m]}
    \includegraphics[width=0.90\columnwidth]{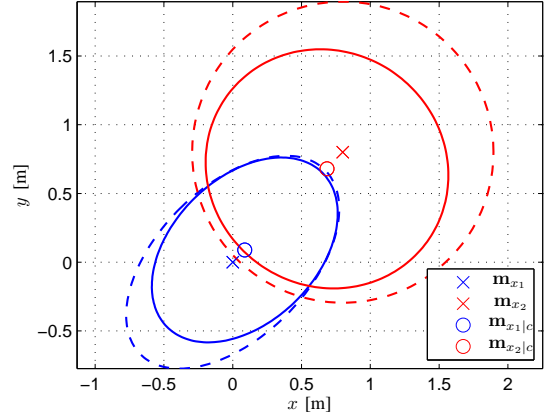}
  \end{center}
  \vspace*{-0.8cm}
  \caption{Illustration of how the side information affects the mean positions. Shown are the means and approximated conditional means, together with the loci corresponding to the 95\% confidence ellipses of Gaussian pdfs with covariance matrices $\mbf{C}_{x_1}$, $\mbf{C}_{x_2}$, $\mbf{C}_{x_1|c}$, and $\mbf{C}_{x_2|c}$, respectively.}
  \label{fig:x}
\end{figure}

\begin{figure}[!t]
  \begin{center}
    \psfrag{mz1}[c][c][0.7][0]{$\mbf{m}_{z_1}$}
  \psfrag{mz1c}[c][c][0.7][0]{$\mbf{m}_{z_1|c}$}
    \psfrag{xaxis}[c][c][0.7][0]{$\Delta x$ [m]}
    \psfrag{yaxis}[c][c][0.7][0]{$\Delta y$ [m]}
    \includegraphics[width=0.90\columnwidth]{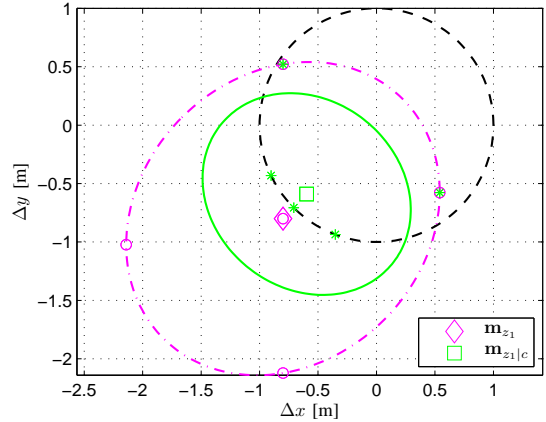}
  \end{center}
  \vspace*{-0.8cm}
  \caption{Illustration of the sampling technique used to approximate $\mbf{m}_{z_1|c}$ and $\mbf{C}_{z_1|c}$.  Sample points before $\mbf{s}^i$ (circles) and after $\mbf{z}^{i}_1$
    (asterisk) the projection, the boundary of the constraint $c$ with $\gamma =1$ (dashed circle), and the confidence ellipses of $\mbf{C}_{z_1}$ (dash-dotted) and $\mbf{C}_{z_1|c}$ (solid) are also shown.}
  \label{fig:z1}
\end{figure}

We illustrate the sigma-point approximation and the resulting
estimator using the following example. Two objects, with a joint spatial Gaussian distribution, are located in $\mathbb{R}^2$. The mean positions of the objects are $\mbf{m}_{x_1}=\mbf{0}_2$ and $\mbf{m}_{x_2}=0.8\cdot\mbf{1}_2$ [m], and the joint state covariance matrix $\mbf{C}_{x}=\mbf{C}_{x_1}\oplus\mbf{C}_{x_2}$, where
\begin{equation*}
\mbf{C}_{x_1}=
\begin{bmatrix}
0.1 & 0.05 \\
0.05 & 0.1
\end{bmatrix} \quad\mbox{and}\quad
\mbf{C}_{x_2}=
\begin{bmatrix}
0.2 & 0 \\
0 & 0.2
\end{bmatrix}.
\end{equation*}
We then provide the side information $c$ with $\gamma = 1$ [m]. Setting $\alpha =
0.95$, the approximated conditional means $\mbf{m}_{x_1|c}$ and $\mbf{m}_{x_2|c}$, along with confidence ellipses, are shown in Fig.~\ref{fig:x}. Fig.~\ref{fig:z1} illustrates the deterministic sampling procedure in subsystem $\mbf{z}_1$. The sample points are generated according to the confidence ellipse and projected orthogonally when the constraint is violated. Observe that the ellipses shrink when the side information $c$ becomes available.

%% file: experiment.tex
\section{Experimental results}
Using Monte Carlo simulations, the estimator $\hat{\mbf{x}}_{|c}$ is evaluated numerically in a positioning and a tracking scenario in $\mathbb{R}^2$. Throughout the experiments, the side information $c$ is given by $\gamma=1$~[m] and the estimator design parameter $\alpha = 0.95$. As a performance measure, the root mean square error (RMSE) of the state estimates is used.

\subsection{Positioning scenario}
\begin{figure}
  \begin{center}
    \psfrag{mx}[c][c][0.7][0]{$\mbf{m}_{x}$}
    \psfrag{xc}[c][c][0.7][0]{$\hat{\mbf{x}}_{|c}$}
    \psfrag{tracesq}[c][c][0.7][0]{$\sqrt{ \text{tr}\{ \mbf{C}_{\hat{x}} \} }$}
    \psfrag{yaxis}[c][c][0.7][0]{RMSE [m]}
    \psfrag{xaxis}[c][c][0.7][0]{$\sigma_1$ [m]}
    \includegraphics[width=0.95\columnwidth]{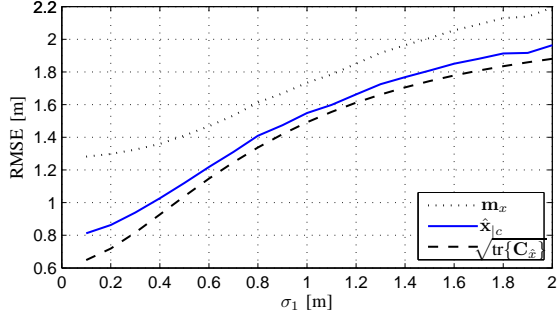}
  \end{center}
  \vspace*{-0.5cm}
  \caption{Positioning scenario. RMSE versus $\sigma_1$ for $\beta = 1$ and side information $c$ with $\gamma=1$.}
  \label{fig:MC_sigma1}
\end{figure}

\begin{figure}
  \begin{center}
    \psfrag{mx}[c][c][0.7][0]{$\mbf{m}_{x}$}
    \psfrag{xc}[c][c][0.7][0]{$\hat{\mbf{x}}_{|c}$}
    \psfrag{tracesq}[c][c][0.7][0]{$\sqrt{ \text{tr}\{ \mbf{C}_{\hat{x}} \} }$}
    \psfrag{yaxis}[c][c][0.7][0]{RMSE [m]}
    \psfrag{xaxis}[c][c][0.7][0]{$\beta$ [m]}
    \includegraphics[width=0.95\columnwidth]{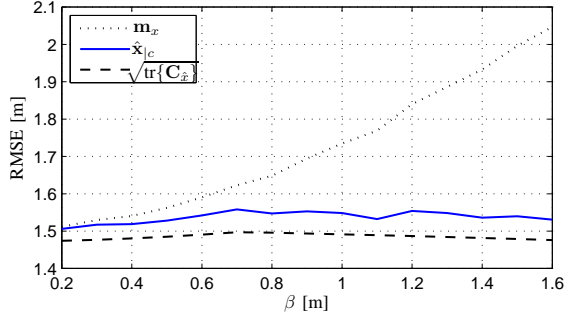}
  \end{center}
  \vspace*{-0.5cm}
  \caption{Positioning scenario. RMSE versus $\beta$ for $\sigma_1 = 1$ and side information $c$ with $\gamma=1$.}
  \label{fig:MC_beta}
\end{figure}

In this scenario, the positioning of the two systems $\mbf{x}_1$ and $\mbf{x}_2$, is considered. In the Monte Carlo simulation, the positions of the two systems were generated by drawing two candidate positions from two independent Gaussian distributions, and then retaining a realization that fulfilled the constraint with $\gamma$. The parameters of the pdf were $\mbf{m}_{x_1} = \beta \mbf{1}_2$, $\mbf{C}_{x_1} = \sigma^2_1 \mbf{I}_2$, and $\mbf{m}_{x_2} = \mbf{0}_2$ and $\mbf{C}_{x_2} = \mbf{I}_2$. The empirical RMSE of the estimator $\hat{\mbf{x}}_{|c}$ was calculated using Monte Carlo simulations with $10^4$ runs. Figs.~\ref{fig:MC_sigma1} and \ref{fig:MC_beta} show the RMSE of the estimator for different values of $\sigma_1$ and $\beta$, respectively. As expected, the accuracy improvement of $\hat{\mbf{x}}_{|c}$ over the prior mean $\mbf{m}_x$ increases as the certainty of the position of one object grows, or when its prior mean is placed further away. In these cases $c$ provides more information.
The figures also show $\sqrt{\mbox{tr}\{\mbf{C}_{\hat{x}} \} }$, which indicate that the second-order statistics are slightly underestimated.

\subsection{Tracking scenario}
Let us now consider the scenario where we want to fuse the position information from two dead-reckoning systems mounted on a non-rigid body, but where the body has a known finite length. This could, for example, be two foot-mounted inertial navigation systems (one system on each foot) used to track the position of a person while walking inside a building~\cite{Brand&Phillips2003}. The joint position state $\mbf{x}(k)= [(\mbf{x}_{1}(k))^\top \quad
(\mbf{x}_{2}(k))^\top]^\top\in\mathbb{R}^4$ of the two dead-reckoning systems is updated every second via the recursion $\mbf{x}(k)=\mbf{x}(k-1) + \mbf{u}(k)$. The displacements $\mbf{u}(k)$ are measured in noise, $\widetilde{\mbf{u}}(k) \sim\mathcal{N}(\mbf{u}(k),\mbf{Q})$, where $\mbf{Q}$ denotes the error covariance of the measured displacements. Due to the integrative nature of the dead-reckoning recursion, the uncertainty of the state $\mbf{x}(k)$ will grow without bound; the covariance of the state uncertainty, $\mbf{P}(k)$, is given by the recursion $\mbf{P}(k)=\mbf{P}(k-1)+\mbf{Q}$, starting from the covariance $\mbf{P}(0)$ of the initial position state.

To reduce the rate at which the position uncertainties grow, we can use the side information that the two systems are mounted on a non-rigid body of known finite length, i.e., we have the side information that $\|\mbf{x}_{1}(k)-\mbf{x}_{2}(k)\|\leq\gamma$, $\forall k$. Given the sequence $\overline{\mbf{u}}(k)=\{\widetilde{\mbf{u}}(i)\}_{i=0}^k$ of measured position changes and the general side information $I$, the conditional pdf $p( \mbf{x}(k) | \overline{\mbf{u}}(k) , I ) $ can be calculated recursively as
\begin{equation*}
\begin{split}
p( \mbf{x}(k) | \overline{\mbf{u}}(k) , I ) =& \int p\left( \mbf{x}(k) | \mbf{x}(k-1), \widetilde{\mbf{u}}(k) \right)\\
&\times p\left( \mbf{x}(k-1) | \overline{\mbf{u}}(k-1) , I \right) d \mbf{x}(k-1),
\end{split}
\end{equation*}
starting from the pdf $p( \mbf{x}_{0} | I )$ of the initial position state. The transition pdf $p( \mbf{x}(k) | \mbf{x}(k-1), \widetilde{\mbf{u}}(k))$ equals $\mathcal{N}( \mbf{x}(k-1) + \widetilde{\mbf{u}}(k), \mbf{Q} )$. The MMSE estimate of $\mbf{x}(k)$ is given by the mean of $p( \mbf{x}(k) | \overline{\mbf{u}}(k) , I )$, which is intractable in general.

When $I = c$, we use the conditional moments at each time instant $k$ and approximate $p( \mbf{x}(k) | \overline{\mbf{u}}(k) , c ) $ by a Gaussian pdf, $\mathcal{N}(\mbf{m}_{x|c}(k), \mbf{C}_{x|c}(k))$. Then the recursion results in Gaussians with computable means that form point estimates $\hat{\mbf{x}}_{|c}(k)$. We compare this to the case when the side information is the actual distances, $I = \{  \mbf{y}(i) \}^{k}_{i=0}$, where $\mbf{y}(k) = \| \mbf{x}_1(k) - \mbf{x}_2(k) \|$.

\begin{figure}
  \begin{center}
    \psfrag{Upper bound}[c][c][0.7][0]{$\sqrt{ \mbox{tr}\{\mbf{P}(k)\} }$}
  \psfrag{Lower bound}[c][c][0.7][0]{$\sqrt{ \mbox{tr}\{\mbf{J}^{-1}(k)\} }$}
  \psfrag{Estimator}[c][c][0.7][0]{$\hat{\mbf{x}}_{|c}(k)$}
    \psfrag{xaxis}[c][c][0.7][0]{Time $k$ [s]}
    \psfrag{yaxis}[c][c][0.7][0]{RMSE [m]}
    \psfrag{title}[c][c][0.7][0]{}
    \includegraphics[width=1.00\columnwidth]{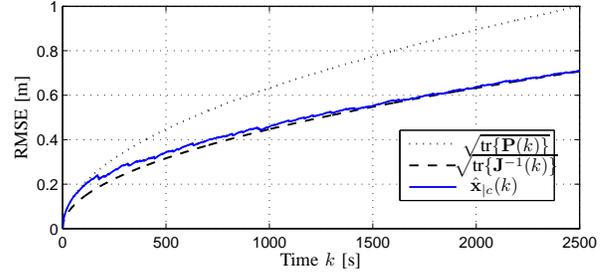}
  \end{center}
  \vspace*{-0.5cm}
  \caption{Tracking scenario. RMSE over time for two dead-reckoning systems when the side information $c$ is used. The PCRBs for the case with no side information and the case with perfect distance information are also shown.}
  \label{fig:filtering}
\end{figure}

The posterior Cram\'{e}r-Rao bound (PCRB) on the RMSE for $I = \varnothing$ is $\sqrt{\text{tr}\{ \mbf{P}(k) \} }$ and for $I = \{  \mbf{y}(i) \}^{k}_{i=0}$, it is $\sqrt{\text{tr}\{ \mbf{J}^{-1}(k)\} }$, where $\mbf{J}(k)$ is the Fisher information matrix of the state (see \cite{TichavskyEtAl1998} for details). The performance of the estimator is compared with the PCRBs in Fig.~\ref{fig:filtering}, where $\mbf{Q} = 10^{-4}\cdot \mbf{I}_4$~[$\text{m}^2$]. Initially, the RMSE follows the upper PCRB, but as the errors of the dead-reckoning systems accumulate, the distance bound becomes more informative and the estimator tends towards the lower PCRB, which has a lower growth rate than the upper PCRB; proof emitted due to space limitations.

\emph{Reproducible research:} The proof of convergence to the lower PCRB  and the Matlab code used in all the simulations is available at \url{www.ee.kth.se/~davez/rr-bayes}.

%% file: conclusions.tex
\section{Conclusion and further work}
We have presented an approximate MMSE estimator that uses a given maximum distance between subvectors as side information. By reducing the dimensionality of the problem, a computationally efficient formulation was given. The estimator has a range of potential applications in positioning and localization.
Our further work includes extending the framework to larger systems with several distance bounds, and applying it to a multi-user indoor navigation system.